\documentclass{article}
 \title{On the Abel-Radon transform of locally residual currents}
\author{Bruno Fabre} 

\begin{document}
\maketitle
   \def \N{ \mbox{I\hspace{-.15em}N}}
   \def \Z{ \mbox{Z\hspace{-.3em}Z}}
   \def \R{\mbox{I\hspace{-.17em}R}}
   \def \P{ \mbox{I\hspace{-.17em}P}}
   \def \C{ \mbox{l\hspace{-.47em}C}}
   \def \Q{ \mbox{l\hspace{-.47em}Q}}

   \newtheorem {thm}{Theorem\hspace{2pt}}
   \newtheorem {prb}{Problem\hspace{2pt}}
   \newtheorem {pro}{Proposition\hspace{2pt} }
   \newtheorem {conj}{Conjecture\hspace{2pt}}
   \newtheorem {cor}{Corollary\hspace{2pt}}
   \newtheorem {defi}{Definition\hspace{2pt}}
   \newtheorem {asum}{Hypothesis\hspace{2pt}}
   \newtheorem {rem}{Remark\hspace{2pt}}
   \newtheorem {lem}{Lemma\hspace{2pt}}
   \newtheorem {exm}{Example\hspace{2pt}}

    \def\CQFD {\hbox{\vrule height 5pt depth 5pt width 5pt}}
       \def\dem{ Proof.}

\begin{abstract}
First we recall the definition of locally residual currents and their basic properties. We prove in this first section a trace theorem, that we use later. Then we define the Abel-Radon transform of a current ${\cal R}(\alpha)$, on a projective variety $X\subset \P^N$, for a family of $p-$cycles of incidence variety $I\subset T\times X$, for which $p_1:I\to T$ is proper and $p_2:I\to X$ is submersive, and a domain $U\subset T$. Then we show the following theorem, for a family of sections of $X$ with $r-$planes (which was proved for the family of lines of $X=\P^N$ by the author in \cite{Fabre1} for $p=1$, in \cite{Fabre2} for ${\cal R}(\alpha)=0$ and $p-$planes for any $q>0$, and by Henkin and Passare in \cite{HenPas} for $p-$planes in $\P^N$ and integration currents $\alpha=\omega\wedge[Y]$, with a meromorphic $q-$form $\omega$, and projective convexity on $\tilde U$):

{\it Let $\alpha$ be a locally residual current of bidegree $(q+p,p)$ on $U^*$, with $U^*:=\cup_{t\in U}{H_t}\subset X$, where ${t}\times H_t:=p_1^{-1}(t)$. Then ${\cal R}(\alpha)$ is a meromorphic $q-$form on $U$, holomorphic iff $\alpha$ is $\overline{\partial}-$closed.

 Let us assume that $\alpha$ is $\overline\partial-$closed, and $q>0$. If ${\cal R}(\alpha)$ extends meromorphically (resp. holomorphically) to a greater domain $\tilde U$, then $\alpha$ extends in a unique way as a locally residual current (resp. $\overline\partial-$closed) to the greater domain ${\tilde U}^*\subset X$.
}

In particular we recover the result of \cite{Fabre2} without using \cite{HenGin}.
We formulate another generalization, for complete intersections with respect to a fixed multidegree.
\end{abstract}

\section{Locally residual currents}

 Let us recall what a locally residual current is on a complex manifold $X$ of dimension $n$ (cf. \cite{Herrera}).
Let us be given $p+1$ complex hypersurfaces $Y_1,\dots,Y_{p+1}$, such that for $i,1\le i\le p+1$, $Y_1\cap \dots\cap Y_i$ is of pure codimension $i$.
Then we define, for a meromorphic $q-$form $\Psi$ with $Pol(\Psi)\subset Y_1\cup\dots\cup Y_{p+1}$ the currents $Res_{Y_1,\dots,Y_{p+1}}(\Psi)$ (resp. $Res_{Y_1,\dots,Y_{p}}P_{Y_{p+1}}(\Psi)$, or simply $P_{Y_1}(\Psi)$ or $[\Psi]$ the {\it principal value} if $p=0$) as follows. On an open subset $U\subset X$, such that every $Y_i\cap U= \{f_i=0\}$ for some holomorphic function $f_i$ on $U$, we define, for a test-form $\phi$:
$$Res_{Y_1,\dots,Y_{p+1}}(\Psi)(\phi)=\lim_{m\to\infty}{\int_{\vline f_1 \vline=\epsilon_{m,1},\dots,\vline f_p \vline=\epsilon_{m,p},\vline f_{p+1} \vline= \epsilon_{m,p+1}}{\phi\wedge\Psi}},$$
and $$Res_{Y_1,\dots,Y_{p}}P_{Y_{p+1}}(\Psi)(\phi)=\lim_{m\to\infty }{\int_{\vline f_1 \vline=\epsilon_{m,1},\dots,\vline f_p \vline=\epsilon_{m,p},\vline f_{p+1} \vline\ge \epsilon_{m,p+1}}{\phi\wedge\Psi}}$$
for any sequence $\epsilon_m=(\epsilon_{m,1},\dots,\epsilon_{m,p+1})\in \R_+^{p+1}$ such that $lim_{m\to\infty}{\epsilon_{m,1}}= 0$ and $lim_{m\to\infty}{\epsilon_{m,i+1}/\epsilon_{m,i}^k}= 0 (1\le i\le p)$ for all integers $k$.
It does not depend on the sequence, nor on the choice of the $f_i$; thus we can define, by a covering of $X$ and a partition of unity subordinated to it, currents on the whole $X$: we call them residual currents.
The first gives a current of bidegree $(q,p+1)$, the second a current of bidegree $(q,p)$, and moreover we have: $Res_{Y_1,\dots,Y_{p+1}}(\Psi)=\overline{\partial} Res_{Y_1,\dots,Y_{p}}P_{Y_{p+1}}(\Psi)$.
A locally residual current is a current which is locally a residual current.

There is some characterization of a locally residual current $\alpha$, by \cite{Bjork}:

1. $\alpha$ is of bidegree $(q,p)$, with an analytic support $Y$ of codimension $p$;

2. $\overline{{\cal I}_Y}\alpha = 0$;

3. $\overline{\partial} \alpha =0$, outside an analytic hypersurface $S$ of $Y$;
and, if this hypersurface $S$ is not empty, the fourth:

4. $\alpha$ is of standard extension through $S$, i.e.
$\alpha(\phi)=lim_{\epsilon\to 0}{\alpha(\xi_\epsilon\phi)}$,
where $\xi_\epsilon$ is a cut-function, $0$ in a neighborhood of $X$ (see \cite{Bjork}).

In particular, integration currents $\omega\wedge[Y]$ are locally residual currents.

Let be $p: X\to Y$ an analytic morphism of analytic varieties. We have:

\begin{lem}
If $\alpha$ is a locally residual current on $X$, of bidegree $(r,s)$ such that $p$ is proper on the support $Z:=Supp(\alpha)$ then $\alpha':=p_*(\alpha)$ is locally residual on $Y$, of bidegree $(r-k,s-k)$ with $k:=dim(Y)-dim(X)$ and support $\subset Z':=p(Z)$, $\overline{\partial}-$ (resp. $\partial-$)closed if $\alpha$ is.
\end{lem}
\dem  First we see that $\alpha'$ is of bidegree $(r-k,s-k)$ by the definition of $p_*$, since the pull-back of forms doesn't change the bidegree. And $\alpha'$ is $\overline{\partial}-$closed if $\alpha$ is, since $p_*$ commutes with $\overline{\partial}$.  We know by Remmert that $Z'$ is an analytic subvariety of $Y$. Moreover, we can see that $\overline{{\cal I}}_{Z'}$ annihilates $\alpha'$, since $p^*({\cal I}_{Z'})\subset {\cal I}_Z$. Moreover, the support of $\overline{\partial} \alpha'$ will be in $Z'$. Finally, $\alpha'$ will remain of standard extension. From this and from the preceding characterization we see that $\alpha'$ is a locally residual if $\alpha$ is.\CQFD

We can see also:
\begin{lem}
If $p$ is a submersion, then: if $\alpha$ is locally residual on $Y$, then $p^*(\alpha)$ is locally residual on $X$, $\overline{\partial}-$ (or $\partial-$) closed if $\alpha$ is.
\end{lem}
\dem It suffices to see it on projections, locally. Then we see by Fubini that: $$p^*(Res_{f_1,\dots,f_q}{\frac{\Psi}{f_1\dots f_q} })=Res_{p^*(f_1),\dots,p^*(f_q)}{\frac{p^*(\Psi)}{p^*(f_1)\dots p^*(f_q)}}$$
Moreover, since $p_*$ commutes with $\overline{\partial}$ for compact supported forms, $p^*$ commutes with $\overline{\partial}$ for currents.
\CQFD

\section{Trace theorem}
Let $\alpha$ be a locally residual current of bidegree $(q+p,p)$ on $D':=D\times \C^p$, with $D\subset \C^n$ proper in the sense that the canonical projection $\pi: D'\to D$ is proper on its support.
Then, by the preceding $\pi_*(\alpha)$ is a $(q,0)-$locally residual current, that is, the principal part of a meromorphic $q-$form $\omega$ which is holomorphic at the places over which $\alpha$ is $\overline{\partial}-$closed.
Let us assume $q=n$. We can for all monomials $m_I\in \C[y_1,\dots, y_p]$ define the traces $u_I$ by $u_I dx_1\wedge\dots\wedge dx_n=\pi_*(\alpha m_I)$, which are meromorphic functions.
Let us denote $Y:=Supp(\alpha)$, and $S:=Supp(\overline{\partial}\alpha)$, $S':=\pi(S)$.
 They can be computed by residues:
 \begin{lem} Outside $S'$, $u_I=\sum_i{Res_{P_i}(\phi_i(x,y)m_I/(f_{1,i}\dots f_{p,i}))}$, where we have put $\alpha$ as an explicit residue at each point $P_i$ of the fiber $\pi^{-1}(x)\cap Y$, with $Y:=Supp(\alpha)$, and we compute it on the fiber $\C^p$ as a punctual residue.
 \end{lem}
 \dem It follows of Fubini formula for the integral and a passage to limit.
 \CQFD

 \begin{lem} Assume that the traces $u_I$ are all zero. Then $\alpha=0$.
 \end{lem}
 \dem We have to use the explicit expression of the trace as a residue: $u_I=\sum_i{Res_{P_i}(\phi_i(x,y)m_I/(f_{1,i}\dots f_{p,i}))}$. Then we see by Fubini that the coefficient at $x$ $\sum_i{Res_{P_i}(\phi_i(x,y)m_I/(f_{1,i}\dots f_{p,i}))}$ is zero for any monomial. But this imply that the current $\alpha_x=\sum_i{Res_{P_i}(\phi_i(x,y)/(f_{1,i}\dots f_{p,i}))}$ is zero on $\C^n$, since we can find for each smooth function $\phi$ an interpolation polynomial $P$ for which $\alpha_x(\phi)=\alpha_x(P)$. So each $\alpha_x$ is zero, and by Fubini $\alpha$ is zero.
 \CQFD

Let us assume that $\alpha$ is $\overline{\partial}-$closed.
Then, we can write explicitly $\alpha$ by using the traces.
For this, let us consider the different projections $\pi_i:\C^{n+p}:(x_1,\dots,x_n,y_1,\dots,y_p)\to (x_1,\dots,x_p,y_i)$. The support $Y$ of $\alpha$ is projected into an hypersurface $Y_i$. There is a minimal polynomial, $P_i=y_i^{d}+a^i_1(x) y_i^{d-1}+\dots+a^i_d(x)=0$, for which $P_i\alpha=0$; $P_i$ is in the ideal of $Y$.
Thus $\alpha,\alpha y_i,\dots, y_i^{d-1}\alpha$ are independant over ${\cal O}(D)$, and $\alpha,\alpha y_i,\dots, y_i^{d}\alpha$ are dependant.
\begin{lem} The $a^i_j$ can be computed from the traces $u_I$.
\end{lem}
\dem First we get the following relations from $y_1^{k_1}\dots y_p^{k_p}P_i\alpha=0$, by taking the trace $\pi_*$:
$$u_{k_1,\dots,k_i+d,\dots,k_p}+a^i_1(x)u_{k_1,\dots,k_i+d-1,\dots,k_p}+\dots+
a^i_d(x)u_{k_1,\dots,k_i,\dots,k_p}=0$$
We see that, since $\alpha,\alpha y_i,\dots, y_i^{d-1}\alpha$ are independant, so are the columns of the matrix $u_{k_1,\dots,k_i+d-1,\dots,k_p},\dots,u_{k_1,\dots,k_p}$, where the lines correspond to multiindices $(k_1,\dots,k_p)$. So in the equations $$a^i_1(x)u_{k_1,\dots,k_i+d-1,\dots,k_p}+\dots+a^i_d(x)u_{k_1,\dots,k_i,\dots,k_p}=-u_{k_1,\dots,k_i+d,\dots,k_p}$$
we can by varying the multi-indices find some non-identically zero determinant, and so solve meromorphically the equations on $(a^i_1,\dots,a^i_d)$, and this for all $i$.
\CQFD

\begin{lem} The current $\alpha$ can be expressed explicitly as a global residue from the traces $u_I$ and the coefficients $a^i_j$.
\end{lem}
\dem Let us consider the following formal series: $\phi=\sum_{k_1,\dots,k_p}{u_{k_1,\dots,k_p}/(y_1^{k_1+1}\dots y_p^{k_p+1})}$. Then we can show by computing that, formally, $P_1\phi$, if we develop, doesn't contain any negative power of $y_1$, by using the relations $$u_{k_1+d,\dots,k_p}+a^1_1(x)u_{k_1+d-1,\dots,k_p}+\dots+a^1_d(x)u_{k_1,\dots,k_p}=0$$
The result is the following: for $k_1<d$ the coefficient of $y_1^{d-1-k_1}$ is:
$(u_{k_1,\dots,k_p}+\dots+a^1_1 u_{k_1-1,\dots,k_p}+\dots+a^1_{k_1} u_{0,k_2,\dots,k_p})/(y_2^{k_2+1}\dots y_p^{k_p+1})$. So we have different formal series $a^1_{k_1} u_{0,k_2,\dots,k_p})/(y_2^{k_2+1}\dots y_p^{k_p+1})$,\dots  $u_{k_1,\dots,k_p}/(y_2^{k_2+1}\dots y_p^{k_p+1})$.
By multiplying by $P_2$ and using the relations $$u_{k_1,k_2+d,\dots,k_p}+a^2_1(x)u_{k_1,k_2+d-1,\dots,k_p}+\dots+a^2_d(x)u_{k_1,\dots,k_p}=0$$
we succeed in eliminating every negative power in $y_2$. So we get that $P_1\dots P_p\phi$ is a polynomial $Q$ in $y_1,\dots,y_p$, with first term $y_1^{d-1}\dots y_p^{d-1}$. So $\phi$ is also, with respect to the $y_i$, a rational function; and so is it with respect to the variables $y_i'=1/y_i$. But with respect to this variables, we get $0$ at $y_i'=0$; so it is well-defined at the origin, and we can develop with respect to the $y_i'$; the series converges, and it has to be the same as the one given by $\phi$. This one is not only a formal series, but thus converges at infinity. Thus we can distribute the integral $\int_{\vline y_1\vline =R_1,\dots,\vline y_p\vline =R_p}$ terms by terms in the sum, seeing we using the variables $y_i'=1/y_i$ that the only non-zero term is $u_{k_1,\dots,k_p}$. But this integral for $Q/(P_1\dots P_p) dy$ is also the same as the sum of the integral $\int_{\vline P_1\vline \epsilon_1,\dots,\vline P_p\vline =\epsilon_p}$, at the different common zeroes of the $P_i$, since we can find between ${\vline y_1\vline =R_1,\vline P_2\vline =\epsilon_2\dots,\vline P_p\vline =\epsilon_p}$ and ${\vline P_1\vline =\epsilon_1,\dots,\vline P_p\vline =\epsilon_p}$ a real variety of dimension $(p+1)$ outside $P_1=0$, over which the form is closed; and repeat it. Thus this is the trace of $\alpha'= Res_{P_1=0,\dots,P_p=0}{Q/(P_1\dots P_p) dy}$. So $\alpha'$ has the same traces as $\alpha$, and so it is equal to $\alpha$ by the preceding.
\CQFD

\begin{thm}
Let us assume that the traces $u_I=\pi_*(\alpha m_I)$, holomorphic on $D$, extend holomorphically to a greater domain $\tilde D$, for all monomials $m_I\in \C[y_1,\dots, y_p]$. Then $\alpha$ extends in a unique way to a $\overline{\partial}-$closed proper locally residual current of bidegree $(n+p,p)$  on ${\tilde D}':={\tilde D}\times \C^p$.
\end{thm}
\dem
First assume that the $u_I$ extend holomorphically. Then we have some determinant over $D$ of the equations:
$$a^i_1(x)u_{k_1,\dots,k_i+d-1,\dots,k_p}+\dots+a^i_d(x)u_{k_1,\dots,k_i,\dots,k_p}=-u_{k_1,\dots,k_i+d,\dots,k_p}$$
which is non-identically zero, and so it is also non-identically zero on $\tilde D$; so the $a^i_j$ extend meromorphically.
We can as before define the function $\phi= Q/(P_1\dots P_p)$ which is rational in the $y_i$. We can develop it, since it is well-defined at infinity, and see that the push-forward of $Res_{P_1=0,\dots,P_p=0}{\phi y_1^{i_1}\dots y_p^{i_p} dx\wedge dy}$ is $u_{k_1,\dots,k_p}dx$. Thus it extends $\alpha$.
The $a^i_j$ are in general meromorphic; but we can multiply them by a common multiple of the denominators to express $\phi= Q/(P_1\dots P_p)$, if the $u_I$ are holomorphic, with holomorphic functions. If the $u_I$ extend meromorphically, the numerator $Q$ is no more holomorphic, and the corresponding current $Res_{P_i=0}{ Q/(P_1\dots P_p)}dx\wedge dy$ is no more $\overline\partial-$closed.
\CQFD

\section{Pull-back, push-forward, and Radon transform}
Let us consider a submersion $p:Y\to X$. The push-forward of a current $\alpha$, on the support of which $p$ is proper, is defined by: $p_*(\alpha)(\phi)=\alpha(p^*(\phi))$. We see that $deg(p_*(\alpha))=deg(\alpha)-p$, where $p=deg(Y)-deg(X)$. Let us assume that $\alpha=[\Psi]$ is associated to a smooth form. Then $p_*([\Psi])$, defined by "integration on the fibers", is also associated to a smooth form, which will be denoted $p_*(\Psi)$. We have $deg(p_*(\Psi))=deg(\Psi)-p$ ($p_*(\Psi)=0$ if $deg(\Psi)<p$). Let us assume $deg(\Psi)=p+q$. Then we can write $p_*(\Psi)=\sum_{i_1<\dots<i_q}{\omega_{i_1,\dots,i_q}dx_{i_1}\wedge\dots\wedge dx_{i_q}}$, where $\omega_{i_1,\dots,i_q}=\int_{p^{-1}(x)}{\Psi/dx_{i_1}\wedge\dots\wedge dx_{i_q}}$.
Here, $\Psi/dx_{i_1}\wedge\dots\wedge dx_{i_q}$ denotes the coefficient of $dx_{i_1}\wedge\dots\wedge dx_{i_q}$ in the decomposition of $\Psi$ with coordinates $x_1,\dots,x_n,y_1,\dots,y_p$.

\section{Abel-Radon transform and generalized Abel's theorem for locally residual currents}

Let us consider a projective variety $X\subset \P^N$, and a family of $p-$cycles $I\subset T\times X$, where we suppose $p_1:I\to T$ proper and $p_2:I\to X$ submersion. For any domain $U\subset T$ we associate the dual domain $U^*:=\cup_{t\in U}{H_t}$, with ${t}\times H_t=p_1^{-1}(t)$. We denote $I_U:=p_1^{-1}(U)$. Then $p_1, p_2$ restrict to $I_U$ and we define the Abel-Radon transform of a current $\alpha$ on $U^*=p_2(I_U)$ by the formula:
${\cal R}(\alpha):=(p_1)_*(p_2^*(\alpha))$, which is a current on $T$.

It follows from the commutation for push-forward and pull-back that the transformation $\cal R$ commutes with $d,\partial, \overline{\partial}$.
Moreover, we can use change of parameters:

Let $\mu: T'\to T$ be change of parameter, such that the new family $I'\subset T'\times X$ with $(t',x)\in I'$ if $(\mu(t),x)in I$, is submersive ($T$ can be, for instance, a submersion, or an appropriate subset of $T$). We can consider the Abel-Radon transform ${\cal R}'$ with respect to this new family. Then we have
\begin{lem}
${\cal R}'(\alpha)=\mu^*({\cal R}(\alpha))$.
\end{lem}
\dem Let us define the map: $\mu':I'\to I$ by $\mu'(t',x)=(\mu(t),x)$. Then we have: $\mu\circ  p_1'=p_1\circ\mu'$. Then we have ${\cal R}'(\alpha)={p_1'}_*(\mu'^*(p_2^*(\alpha)))$, so that we have to show $\mu^*\circ {p_1}_*={p_1'}_*\circ\mu'^*$. If we compose with $\mu'_*$, we have to check by the first relation that $\mu^*\circ\mu_*\circ {p_1}_*=p_1^*\circ\mu'^*\circ\mu'_*$.
\CQFD

We can also change the variety $X$: Let $\mu: X'\to X$ be an injection, and the corresponding family $(t,x)\in I'$ if $(t, \mu(x))\in I$, so that $I'$ is submersive on $I'$, and of transformation ${\cal R}'$. Then:
\begin{lem}
${\cal R}'(\alpha)={\cal R}(\mu_*(\alpha))$.
\end{lem}
\dem ${p_1}_*(p_2^*(\alpha))={p_1'}_*\circ mu'_*\circ p_2^*\circ \mu^*\circ\mu_*(\alpha)$ by the isomorphism, which is by the commuting diagram ${p_1'}_*(p_2'^*(\mu_*(\alpha)))$
\CQFD

\subsection{An explicit expression of the transform in an affine chart}
Let us assume that the form $\alpha$ is of bidegree $(p,p)$.
Then: ${\cal R}(\alpha):=\int_{H_t}{\alpha}$. We can see this by the fact that:
${\cal R}(\alpha)=\int_{p_1^{-1}(t)}{p_2^*(\alpha)}$; but the integral of $p_2^*(\alpha)$ over ${t}\times H_t$ is nothing but $\int_{H_t}{\alpha}$.
Let us now assume that $\alpha$ is of bidegree $(n+p,p)$.
Then:
 \begin{lem}${\cal R}(\alpha)=\sum_{0\le i_j\le p}{u_{i_1,\dots,i_n}da^1_{i_1}\wedge\dots\wedge da^n_{i_n}}$, where:
$u_{i_1,\dots,i_n}=u^{j_0,\dots,j_p}=\int_{H_t}{Res_{l_1=0,\dots,l_p=0}{\alpha Y_0^{j_0}Y_0^{j_0}\dots Y_p^{j_p}}/l_1\dots l_n}$, $j_i:=card\{k,i_k=i\}$, $j_0+\dots+j_p=n$.
Here, we use the following affine coordinates:
$X_1=\sum_{j=0}^p{a^1_j Y^j},\dots, X_n=\sum_{j=0}^p{a^n_j Y^j}$; we let $a^i_0=b_i$.
\end{lem}
\dem First let us notice that the residue is a well-defined form of bidegree $(p,p)$ on $H_t$.
We have defined the ${\cal R}(\alpha)$ as ${p_1}_*(p_2^*(\alpha))$; but this is also equal to ${\pi_1}_*(\pi_2^*(\alpha)\wedge [I])$, where $I$ is the incidence variety, defined by the linear equations $l_1=0,\dots,l_n=0$, and $\pi_1:T\times X\to T, \pi_2:T\times X\to X$ are the canonical projections. But this also equal to ${\pi_1}_*(Res_{l_1=0,\dots,l_n=0}{\pi_2^*(\alpha)/{l_1\dots l_n}}\wedge dl_1\wedge\dots\wedge dl_n)$. Let us write $dl_1\wedge\dots\wedge dl_n$ in the form $\sum_{0\le i_j\le p}{Y_0^{j_0}\dots Y_p^{j_p} da^1_{i_1}\wedge\dots\wedge da^n_{i_n}+m}$, with $j_0+ \dots+j_p=n$, and $m$ contains differentials with $X_1,\dots,X_n,Y_0,\dots,Y_p$. Then the term with $m$ will vanish in $\pi_2^*(\alpha)\wedge dl_1\dots dl_n$. It will stay $\sum_{0\i_j\le p}{\int_{H_t}{Res_{l_1,\dots,l_n}{\alpha/{l_1\dots l_n}Y_0^{j_0}\dots Y_p^{j_p}}}da^1_{i_1}\wedge\dots\wedge da^n_{i_n}}$, which is what we wanted.
\CQFD

When $\alpha$ is a locally residual current, we can express the coefficient $u_{j_0,\dots,j_p}$ in a different way: it is equal to:
$$ \sum_{i}{Res_{P_i}{\Psi Y_0^{j_0}\dots Y_p^{j_p}/{Q_1\dots Q_p l_1\dots l_n}}}$$
where the $P_i$ are the intersections of $H_t$ with the support of $\alpha$, and we write at $P_i$ $\alpha=Res_{Q_1,\dots,Q_p}{\Psi/{Q_1\dots Q_p}}$. This punctual residue doesn't depend on the way we express $\alpha$ as a residue.

\section{The inverse Abel's theorem}

We consider as family $T:=G(p, N)$ the grassmannian of $p-$planes in $\P^N$, with $I$ the incidence variety, with projections $p_1:I\to T, p_2:I\to \P^N$, and $U\subset T$ some domain, with dual $U^*$. If $X\subset \P^N$ is a projective variety, we have by the preceding: ${\cal R}(\alpha):= {p_2}_*(p_1^*(i_*(\alpha)))$, where $i: V^*\to U^*$ is the inclusion of $V^*=U^*\cap X$ in $U^*$.

Assume now that $\alpha$ is a locally residual current of bidegree $(q+r,r)$ on $V^*$, where $r$ is the generic dimension of $H_p\cap X$. Then we obtain that ${\cal R}(\alpha)$ is the principal value of a meromorphic $q-$form: it is a locally residual current of bidegree $(q,0)$. Then we have:
\begin{thm}
If ${\cal R}(\alpha)$ extends holomorphically (resp. meromorphically) to a greater domain $\tilde U$, then $\alpha$ extends in a unique way as a locally residual current to the greater domain ${\tilde U}^*\cap X$.
\end{thm}

\subsection{Restriction to the case of $X=\P^N$}
We consider the current $\alpha':=i_*(\alpha)$, where $i:V^*\to U^*$ is the natural inclusion. Then we extends $\alpha'$ to $\tilde{\alpha'}$ on ${\tilde U}^*$ if we know that the theorem is valid on $X=\P^N$. But it is clear that the support of $\alpha'$ must remain in $X$: if $f$ annihilates on $X$, $f{\tilde{\alpha'}}=0$ on $U^*$, and thus on ${\tilde U}^*$. Thus we can restrict to the case of $X=\P^N$.

\subsection{Restriction to the case of maximal degree $q=n$}
Let us assume that we have shown the theorem for $q=n$.
Let us assume $\alpha$ of bidegree $(q+p,p)$, with $0<q<n$, where $n:=dim(Y),Y:=Supp(Y)\subset U^*$.
We fix $k:=n-p$ hyperplanes, so by restriction we get a current $\alpha'$ of same bidegree $(q+p,p)$ in some $H=\P^{q+p}$. The Abel-Radon transform of $\alpha'$ with respect to the $p-$planes contained in $H$ correspond to the Abel-Radon transform of $\alpha$ restricted to the corresponding Schubert cycle in $G(p,N=n+p)$. Thus we can prolongate any restriction. If we consider the support $Y$ of $\alpha$, we see that the sections $Y\cap H_t\subset U$ extend to $U^*$. Since this is true for any section with such $H_t$, we see that $Y$ extends to $U^*$, by a variant of the theorem of Harvey-Lawson on the boundaries of projective varieties. Finally, we see that $\alpha$ extends locally, if we write it as a global residue $\sum_I{\Psi_I(x)dx^I}/{f_1\dots f_p}$ where the $f_i$ extends; then we see that, for different values of $\alpha^I$, $\sum_I{\alpha_I \Psi_I}$ extend on $(q+p)-$planes. By taking different $(p+q)-$planes through a point, we can find the $\Psi_I$. Thus $\alpha$ extends locally: at a neighborhood point, we can write: $\alpha=Res_{f_1,\dots, f_p}{\frac{\sum_I{\Psi_I dx^I} }{f_1\dots f_p} }$.
 Then we can consider a maximal prolongation domain $D\subset {\tilde U}$; and this domain must contain $\tilde U$.

\subsection{The case of maximal degree }

Assume now $q=n$. We choose affine coordinates, associated with an hyperplane at infinity; with these coordinates, we can write a $p-$plane in the form:
$l_i=x_i-\sum_{j=1}^p{a_i^j y_j}-b_i=0 (1\le i\le n)$, so we get affine coordinates $a_i^j,b_i$ on the grassmannian $G(p,N)$.
We denote $u_{i_1,\dots,i_p}=\sum_{i}{Res_{P_i}{\Psi y_1^{i_1}\dots y_p^{i_p}/Q_1\dots Q_p l_1\dots l_n}}$; where $\alpha=Res_{Q_1\dots Q_p}{\Psi/{Q_1\dots Q_p}}$. Then we can see that the coefficients of ${\cal R}(\alpha)$ can be written in the form $u_I$: let us look at the explicit expression of the Abel-Radon transform on $P^*$ in affine coordinates. We get $$\sum_{0\le j_1,\dots,j_n\le p}{ {u_I} da_1^{j_1}\wedge\dots\wedge da_n^{j_n}}$$ where $u_I=\sum_i{Res_{P_i}{y^I \frac{\alpha}{l_1\dots l_n}}}$, $I=(i_1,\dots, i_p)$, $a_j^{p+1}:=b_j$ and $y^I=y_1^{i_1}\dots y_p^{i_p}=\frac{D(l_1,\dots, l_n)}{D(a_1^{j_1}\dots a_n^{j_n})}$. Here the $P_i$ are the intersection points of the support of $\alpha$ with the $p-$plane $H_{a,b}$. The residue $Res_{P_i}{\frac{\alpha}{l_1\dots l_n} }$ is by definition $Res_{P_i}{\Phi/{Q_1\dots Q_p l_1\dots l_n}}$, if we write at $P_i$ $\alpha=Res_{Q_1\dots Q_p}{\Phi/{Q_1\dots Q_p}}$.
Then:
\begin{lem}
Let us assume that $\alpha$ is defined and $\overline\partial-$closed on $P^*$, where $P=P_a\times P_b$ is a polydisc; we assume moreover that on $P^*$, the support of $\alpha$ doesn't meet the hyperplane at infinity. If $u_0$ extends holomorphically to a greater polydisc $P'=P_a\times P_b'$, then for all multiindices $I=(i_1,\dots,i_p), i_j\ge 0$, $u_I$ extends holomorphically also to $P'$.
\end{lem}
\dem
We have that ${\cal R}(y^I\alpha)$ is closed for all $I$ on $P$, so we have the following relations in $P$: $\partial_{b_i}{u_{i_1+1,i_2,\dots,i_p}}=\partial_{a_{i}^1}{u_I }$ for all $1\le i\le n$. Since ${\cal R}{\alpha}$ extends on $P'$, we have that $u_{0,\dots,0}$, the coefficient of $b_1\wedge\dots\wedge b_n$, extends on $P'$, and also all $u_{i_1,\dots,i_p}$ for $i_1+\dots+i_p\le n$.
 So by integrating on $b_i$ in $P'$, for fixed $a$ in $P_a$, $u_{i_1+1,i_2,\dots,i_p}$, which is defined on $P_b$, extends also holomorphically in $b$. In fact, $ d u_{i_1+1,i_2,\dots,i_p}=\sum_{i=1}^n{\partial_{a_{i}^1}{u_I }db_i}$ extend on $P_b$, and so since $P_b$ is simply connected, $u_{i_1+1,i_2,\dots,i_p}$ can be extended, if $u_I$ can be extended.
But we know that if a function, holomorphic in $P=P_a\times P_b$, extends holomorphically for $a$ fixed to $P'=P_a\times P_b'$, then it extends holomorphically to $P'$.
So we have shown that if $u_{i_1,\dots,i_p}$ can be extended, $u_{i_1+1,\dots,i_p}$ can also be extended. We could show in the same way that $u_{i_1,i_2+1,\dots,i_p}$ could be extended.
By iteration, we have that if $u_0$ extends to $P'$, $u_I$ extends holomorphically to $P'$ for all $I=(i_1,\dots,i_p)$.
\CQFD

Then by the trace theorem we have that $p_2^*(\alpha)$ extends over $P'$, outside infinity. But on $I$, with coordinates $(x,y,a)$, if we fix $a$, extending $p_2^*(\alpha)$ when we fix $a$ is equivalent to extend $\alpha$ on the corresponding domain. Let us fix an hyperplane at infinity $H_0$.
We can write $p_2:(x,y,a)\to (x,y)$ outside infinity, so that $p_2^*(\alpha)$ is in fact independent of $a$ in these coordinates; so we get a prolongation of $\alpha$ in the corresponding domain: $\cup_{(a',b)\in P'}{H_{a',b}}\backslash H_{a'}$, where $H_{a'}$ corresponds to the $(p-1)-$plane in infinity, center of projection, associated with the fixation of $a'$ (we check that the dimension of the parameter $a'$ is $np$, the dimension of the grassmannian $G(p-1,N-1)$).
 Since we are in a local reasoning, we can assume that the hyperplane at infinity doesn't meet $Y$, the support of $\alpha$. So by the preceding, we have extended $\alpha$ in a greater domain, corresponding to the union of $U^*$ and the $p-$planes with parameter $a'$.
 But then we can also extend ${\cal R}(\alpha)$ and all the traces to the corresponding domain, union of $p-$planes meeting every component of the prolongation of $\alpha$. So we made a local prolongation, and this prolongation can be made at any boundary point of $U^*$, by a convenient choice of the $(p-1)$plane at infinity. Let us denote that, since $p_2^*(\alpha)$ doesn't depend on $a$, the extension doesn't depend of the choice of the $(p-1)-$plane at infinity.

 Let us consider a maximal domain $D\subset {\tilde U}$ for which $\alpha$ extends in the corresponding domain $D^*$; let us denote $\tilde \alpha$ this extension. If $D\not= {\tilde U}$ we will construct a greater domain for which $\alpha$ extends: let us take a point $P$ in the boundary of $D$. Then $\tilde\alpha$ extends, by the preceding local prolongation, for a neighborhood of $P$, by a convenient choice of the $(p-1)-$plane at infinity; that is, $\tilde\alpha$ extends in the corresponding neighborhood of $H_P$, the $p-$plane of $P$; so we have constructed a greater domain where $\tilde\alpha$ extends, which is impossible by assumption. Thus $D=\tilde U$.

Let us denote ${\tilde U}_{a_0}$ the section of ${\tilde U}$ by $a=a_0$.
 Then we have extended $p_2^*(\alpha)$ over ${\tilde U}_{a_0}$, on $I_{{\tilde U}_{a_0}}=p_1^{-1}({\tilde U}_{a_0})$. If we fix the coefficients $a_{i}^j$, we can determine a section $s:{\tilde U}_{a_0}^*\backslash H\to I_{{\tilde U}_{a_0}}$, where $H={Y_0=0}$ is the hyperplane at infinity, in the following way: to a point $(x,y)$, we associate the $p-$plane defined by $b_1=x_1-a^1_1 y_1-\dots-a^1_p y_p,\dots, b_n=x_1-a^n_1 y_1-\dots-a^n_p y_p$, where $a_0=(a^i_j)$. Then we can define $\tilde \alpha:=s^*(\alpha')$, where $\alpha'$ is the extension of $p_2^*(\alpha)$ on $I_{{\tilde U}_{a_0}}$; it will be an extension of $\alpha$ on ${\tilde U}_{a_0}^*\backslash H$, and moreover it doesn't depend of the choice of the affine chart. If we vary $a_0$ in a polydisc $D$, we get an extension of $\alpha$ on the reunion $\cup_{a_0\in D}{{\tilde U}_{a_0}^*\backslash H}$. Thus we can extend $\alpha$ on the reunion of the cones of ${\tilde U}^*$ associated to a $(p-1)-$plane contained in $U^*$.
By iteration, we can extend $\alpha$ to ${\tilde U}^*$. 

 \subsection{The meromorphic case}
 Let us assume the prolongation of ${\cal R}(\alpha)$ to $\tilde U$ is meromorphic.
 It us sufficient to extend $\alpha$ locally, the general case is a consequence of this. So let us assume that we have extended $\alpha$ for a maximal domain $D$ (which is in particular linearly convex in the sense that it corresponds to all $p-$planes contained in $D^*$). Then let us choose in $D^*$ a boundary point $P$, with a $p-$plane through it. We can choose an hyperplane at infinity outside $P$, so we make "projections" from a $\P^{p-1}$ outside it. But fixing the $\P^{p-1}$ at infinity corresponds to fix $a=(a_i^j)$ in the affine coordinates.
 Let us assume that $D$ doesn't contain ${\tilde U}$. Then we have that some traces $u_I$ which extends in a neighborhood of $P$, in fact the coefficients of ${\cal R}(\alpha)$. But let us consider ${\cal R}(\alpha y_i)$. It is closed, so we get the differential equations: $\partial_{b_i}{u_{i_1+1,i_2,\dots,i_p}}=\partial_{a_{i}^1}{u_I }$ for all $1\le i\le n$. But we know (cf. \cite{Fabre1}) that if $f$ is holomorphic outside $S$, and $df$ is meromorphic through $S$, then also $f$ extend meromorphically through $S$. Thus, for $a$ fixed, $\partial_{b_i}{u_{i_1+1,i_2,\dots,i_p}}$ extends meromorphically through the poles. By induction, all traces $u_I$ extend meromorphically. Then we can use the trace theorem for $a$ fixed to extend $\alpha$ in a neighborhood of $P$, and we get a greater domain than $D$. Thus we must have that $D$ contains ${\tilde U}$.

\section{Applications}
First we get, if $\alpha$ is an integration current $\omega\wedge[Y]$,
and if ${\cal R}(\alpha)$ extends holomorphically, an extension of $\alpha$ which is also an integration current, so we recover the theorem of Henkin and Passare (\cite{HenPas}) in the holomorphic case.
If ${\cal R}(\alpha)=0$, and $\alpha$ locally residual of bidegree $(p+q,p), q>0$, we get an algebraic prolongation of $Y=Supp(\alpha)$, so we get a generalization of the inverse Abel's theorem of Griffiths (\cite{Griffiths}) to the non-reduced case.
Let us notice that this imply that $\alpha$ can be written as a global residue, so $\alpha=\overline\partial \beta$, with $\beta$ locally residual.
Let us assume that $\beta$ is a locally residual current of bidegree $(q+p-1,p-1)$ in a linearly $p-$concave domain, union of $p-$planes. Then the Abel-Radon transform with respect to $(p-1)-$planes gives a meromorphic $q-$form in a domain of the grassmannian $G(p-1,N)$, which contains a $\P^p$ (the $(p-1)-$planes contained in a $\P^p$), and so is concave in the sense of Andreotti; thus, $q-$form extends by \cite{Dingoyan} in a rational form; and by the preceding, $\beta$ also extends on $\P^N$ as an algebraic current.

Let us recall that the following theorem is a consequence of the theorem of Griffiths:

{\it If in a linear web in $\C^n$ (with hyperplane foliations) we have an abelian relation, this web is algebraic (that is, is defined by the hyperplanes cutting points of an algebraic curve in $\P^n$).}
Let us consider a locally residual current $\alpha$ of maximal degree in a projectively convex domain $U^*$, the web defined by the hyperplanes going through a fixed point in $U$, domain of the dual projective space. Then we can see that ${\cal R}(\alpha)$ is the sum of $d$ holomorphic $1-$forms, each corresponding to the sheets of the web. It would be interesting to see how these $1-$forms express with respect of the defining functions $u_i (1\le i\le d)$ of the sheets.

\section{Generalization: the transform with respect to complete intersections}

 Analogous results could be proved by substituting grassmannian by a parameter space of complete intersection of a given multidegree. The theorem was proved in my thesis for integration currents $\omega\wedge [Y]$.
 The principle of the proof is still valid here:
 restriction to the maximal degree by restriction, and then restriction to the case of curves by taking residues on fixed algebraic hypersurfaces of the complete intersection, and letting free just a polynomial. Then, with just one hypersurface variation, we can return to the case of hyperplanes by the Veronese mapping, using the change of variety in the Abel-Radon transform.

\end{document}